\documentstyle[amscd,amssymb,verbatim,12pt]{amsart}

\input xypic
\xyoption{all}

\def\hcorrection#1{\advance\hoffset by #1 }
\def\vcorrection#1{\advance\voffset by #1 }

\vcorrection{-.5in}
\hcorrection{-.5in}

\newcommand{\C}[1]{{\cal#1}} % Calligraphic
\newtheorem{theorem}{Theorem}[section]

\theoremstyle{definition}
\newtheorem{defin}{Definition}[section]

\theoremstyle{definition}

\theoremstyle{remark}
\newtheorem{rem}{Remark}
\newtheorem{notation}{Notation}[section]

\newcommand{\To}{\longrightarrow}

\begin{document}
\pagestyle{plain}
\addtolength{\footskip}{.3in}

%Topmatter
\title{On Categorification}
\author{Lucian M. Ionescu}
\address{Mathematics Department\\Kansas State University\\
             Manhattan, Kansas 66502}
\email{luciani@@math.ksu.edu}
%\thanks{...}
\date{April 1999}

% End topmatter

%=======================================================================
%                      		Abstract
%=======================================================================
\begin{abstract}
We review several known categorification procedures, and introduce
a functorial categorification of group extensions with applications to
non-abelian group cohomology. Categorification of acyclic models and
of topological spaces are briefly mentioned.
\end{abstract}

\maketitle
\tableofcontents

% *************************  end Header  *********************************

%=================================================================
%                   Introduction
%=================================================================
\section{Introduction}\label{S:intro}
The term {\em categorification} was invented by
L. Crane \cite{CF,C} to denote a process of associating
{\em category-theoretic} concepts to {\em set-theoretic}
notions and relations (see also \cite{Baez}).

An often used correspondence is the following:
$$\begin{array}{|c|c|}
\hline
set\ elements & objects\\
\hline
equality\ between & isomorphisms \\
set\ elements & between\ objects\\
\hline
objects & categories\\
\hline
morphisms & functors\\
\hline
equality\ between & natural\ isomorphisms\\
morphisms & of\ functors\\
\hline
\end{array}$$
The role of such a procedure is to introduce new structures,
for example Hopf categories for constructing 4-D TQFTs \cite{CF}.
It provides examples,
for instance monoidal categories with prescribed fusion rules \cite{Dav1},
and and toy models, for example invariants of 3-D manifolds associated to
finite gauge groups as a special case of more general
constructions (\cite{CY}, pp. 12).

We consider categorification as a "bridge" allowing to
transport constructions from {\em set-algebra}
to {\em higher-dimensional algebra} \cite{Baez1},
and allowing to exploit the category theory results for a better
understanding of the former,
e.g. defining the cohomology of monoidal categories \cite{Sa,Dav2}
and then using it to understand non-abelian group cohomology \cite{Io}.

It is also a {\em method} to "force" the right approach in
solving a problem:
$$
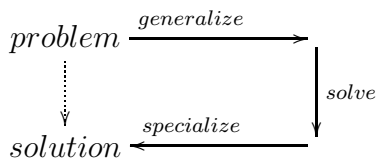
\diagram
problem \ar@{.>}[d] \rrto^{generalize} & \quad & \dto^{solve} \\
solution & \quad & \llto_{specialize}
\enddiagram$$
since categorification may lead to surprising clarifications.

We will consider correspondences as defined below.
\begin{defin}
A {\em categorification} is a functor $\C{C}$ defined on a
concrete category $\C{C}_1$ to a category $\C{C}_2$
with small categories as objects (0-arrows), functors as morphisms (1-arrows),
associating categories to objects and functors to morphisms.
\end{defin}
Often the target is a 2-category with natural transformations as
2-arrows.

The categorification procedures considered may be extended to group rings,
which  provide much more structure as Hopf rings.

%=================================================================
%                   Groups as symmetries 
%=================================================================
\section{Groups as symmetries}
From the beginning of category theory, groups were interpreted as categories
by Eilenberg-Mac Lane \cite{EML} (1942).

\subsection{Tautological categorification}

Let $G$ be a group. Define the category $\C{C}_G$ with one object
$*_G$ having as morphisms the elements of the group.
The group multiplication is used to define the composition of morphisms
in the category. Note that $\C{C}_G$ is a groupoid.
%
% picture explaining briefly the procedure
%
A group homomorphism $f:G\to H$ defines in a tautological way
a functor $\C{C}_f$ between the one object groupoids
associated to the source and target of $f$.
$$\C{C}_f(*_G)=*_H, \qquad \C{C}_f(x)=f(x), \ x\in G$$
The above correspondence is functorial, with domain the category of
groups $\C{G}rp$ and valued in the category $\C{G}rpoid$ of groupoids.
\begin{defin}
The categorification functor $\C{C}$ is called the
{\em tautological categorification}.
\end{defin}
It is an isomorphism between the
categories of groups and one object groupoids.

Natural transformations are one benefit of categorification.
Let $f,g: G\to H$ be group homomorphisms, or in categorical language,
functors between the corresponding one object groupoids.
A natural transformation $\eta:\C{C}_f\to \C{C}_g$ (functorial morphism)
is determined by one element $h=\eta_{*_G}$ of $H$,
satisfying the relation $g(x)=h f(x) h^{-1}$:
$$\diagram
f(*_G) \dto_{f(x)} \rto^{h} & g(*_G) \dto^{g(x)}\\
f(*_G) \rto^{h} & g(*_G)
\enddiagram$$
Then the set of conjugacy classes of group homomorphisms is
the set of equivalence class of functors modulo natural isomorphisms
and it will be denoted by $[\C{C}_G, \C{C}_H]=Hom(\C{C}_G, \C{C}_H)/\sim$.

\begin{notation}
We will call natural isomorphisms of functors {\em homotopies}
and isomorphism classes of functors, homotopy classes.

The map sending an object to its isomorphism class is denoted by $\pi_0$.
The functor associating the homotopy class of a functor
will be denoted by $\pi_1$.
\end{notation}
The relevance of the homotopy class needs not be explained
(fundamental groupoid of a topological space,
holonomy groups, monodromy, etc.).

%
% *************** Groups as abstract symmetries ************************
%
\subsection{Group homomorphisms as representations}

Another benefit of tautological categorification is
the representation of a group as the group of symmetries of an
abstract object. Identifying $G$ with $Aut(*_G)$,
the elements of $G$ are ``abstract'' symmetries of $*_G$.
Here ``identifying'' means the ``identity'' morphism:
$$ G \overset{id}{\To} Aut_{\C{C}_G}(*_G)$$
i.e. the fundamental representation of $Aut(*_G)$ on $*_G$.

More general, the set of group homomorphisms $Hom(G,H)$ can be
interpreted as the category $G-\C{C}_H$ of representations of
the group $G$ on $*_H$, i.e. on the one-object of the category $\C{C}_H$.
Similar to the case of representations of a group in the category of vector
spaces, $\C{V}ect$, the morphisms are $G$-equivariant $\C{C}_H$-morphisms,
 i.e. elements of the group $H$ which intertwine the two morphisms.
In this way conjugacy classes of group homomorphisms are interpreted
as homotopy classes of functors $[\C{C}_G,\C{C}_H]$ or
isomorphism classes of representations.

As an example, we mention the set of {\em outer morphisms}
of the group $H$:
$$Out(H)=Aut(H)/\sim=[H,H]$$ 
and {\em outer actions} $Hom(G, Aut(H))/\sim$ of $G$ on $H$,
as a homotopy class of functors $Hom(\C{C}_G, G-\C{C}_H)$. 

To accommodate quasi-actions and to obtain examples of
monoidal categories, one may interpret the group composition
as a monoidal product, rather then a composition of morphisms.

%=================================================================
%                   Groups as monoidal categories
%=================================================================
\section{Groups as monoidal categories}

The notion of {\em monoidal category} $(\C{M},\otimes,\alpha)$
is itself a categorification. It is the category analog
of a {\em monoid}.

Recall that an equivalence relation
on a set $X$ defines a groupoid with objects the elements of the set
and with a morphism $x\to y$ for each ordered pair of
equivalent elements $x\sim y$.
The composition of morphisms is actually the transitivity property.

There are two important, ``antipodal'' cases \cite{Br}.
\begin{defin}
The {\em discrete groupoid} on the set $X$ is the groupoid
associated to the minimal equivalence relation defined by
the diagonal of $X$. It is denoted by $\C{D}_X$.

The {\em simplicial groupoid} 
on the set $X$ is the groupoid
associated to the maximal equivalence relation defined by
the product $X\times X$. It is denoted by $\C{S}_X$.
\end{defin}

%
% *************** Discrete categorification ************************
%
\subsection{Discrete categorification}
Let $G$ be a group. Define $\C{D}_G$ as the discrete groupoid on the
set $G$. It is a skeletal category, i.e. with only one object in each
isomorphism class. 

There is a natural monoidal product on the category $\C{D}_G$
defined by the multiplication law in the group $G$.
Then $(\C{D}_G, \otimes)$ is a strict monoidal category.

Group homomorphisms define in an obvious way strict monoidal functors.
The correspondence $G \mapsto \C{D}_G$ is functorial, faithful and full,
since in this version of categorification there are no
natural transformations, except identities.
\begin{defin}
The functor $\C{D}:\C{G}rp \to \otimes\text{-}\C{G}rpoid$ is called
{\em discrete categorification}, where $\otimes\text{-}\C{G}rpoid$
is the category of monoidal groupoids.
\end{defin}
It is a construction used in relation to {\em fusion rings},
as it will be briefly mentioned bellow.

%
% *************** K-categorification ************************
%
% maybe expand later
\subsection{K-categorification}

A {\em fusion rule} \cite{Fr,Dav1} is a set $S$ with a
family of non-negative coefficients defining a semi-ring structure
on the free abelian group generated by the set $S$:
the {\em enveloping semi-ring}.

The ``linear version'' of discrete categorification was introduced
by L. Crane and D. Yetter \cite{CY} to provide examples of
semi-simple monoidal categories with a prescribed set of fusion rules.
A {\em K-categorification} of a fusion rule \cite{CY}, with $K$ a field,
is a $K$-linear tensor category with  Grothendieck ring the enveloping
ring corresponding to the fusion rule.
It is a section for the isomorphism class functor
$\pi_0:\otimes\text{-}\C{G}rpoid \to \C{M}on$,
which associates to a monoidal category the monoid of isomorphism classes.

Changing coefficients in the discrete monoidal category associated to a group
by ``tensoring'' with $K$ is a categorification of the group ring $KG$. 
It is a $K$-categorification of the fusion rule defined by
the group ring. The set of such categorifications is in bijection 
with the set of 3-cocycle of the group $G$ with coefficients
in the multiplicative group of the field, as associators \cite{CY}.

%
% *************** Simplicial categorification ************************
%
\subsection{Simplicial categorification}

Let $G$ be a group. Consider the simplicial groupoid $\C{S}_G$.
There is only one isomorphism class, and in fact it is a
{\em contractible groupoid}:
$$
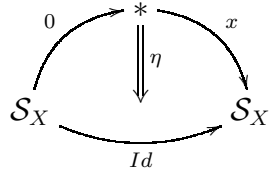
\diagram
& * \ar@{=>}[d]^{\eta} \ar@/^1pc/[dr]^{x} & \\
\C{S}_X \ar@/_1pc/[rr]_{Id}  \ar@/^1pc/[ur]^{0} & & \C{S}_X
\enddiagram$$
i.e. the functor collapsing the category to a chosen base point
is homotopic to the identity functor.

Define the monoidal product on objects as group multiplication.
Since between each pair of objects there is precisely one
morphism only, the extension of the product on morphisms is unique.
We obtain a strict monoidal category $(\C{S}_G, \otimes)$ where
any diagram commutes!

Group homomorphisms define in an obvious way strict monoidal functors
between the corresponding simplicial groupoids,
and the correspondence $G\mapsto \C{S}_G$ is functorial.
\begin{defin}
The functor $\C{S}_G:\C{G}rp\to \otimes\text{-}\C{G}rpoid$ is called
{\em simplicial categorification}.
\end{defin}

%
% *************** Covering transformation ************************
%
\subsection{The covering transformation}\label{covtr}
Define the functor $\C{R}_G:\C{S}_G\to \C{C}_G$ between the simplicial and
tautological categorification of the group $G$:
$$ \C{R}_G(x)=*_G,\quad
\C{R}_G(x\to y)= y x^{-1}, \qquad x, y \in Ob(\C{C}_G)$$
Disregarding the monoidal structures on $\C{S}_G$, the correspondence
defines a natural transformation $\C{R}:\C{S} \to \C{C}$ between
simplicial and tautological categorification functors.

The group $G$ acts freely on $\C{S}_G$ through left and right
multiplication on objects and morphisms.
It also acts on $\C{C}_G$ as inner conjugation from the left and
right trivial.

The functor $\C{R}_G$ is $G$-biequivariant:
{ %\large
$$
\C{R}_G(I_a \otimes \overset{x}{\underset{y}{\downarrow}})=
a\C{R}_G(\overset{x}{\underset{y}{\downarrow}})a^{-1}
$$
$$\C{R}_G(\overset{x}{\underset{y}{\downarrow}} \otimes I_a)=
\C{R}_G(\overset{x}{\underset{y}{\downarrow}})
$$
}
while the corresponding relation on objects is trivial.

%=================================================================
%                   Applications to group extensions
%=================================================================
\section{Applications to group extensions}

Recall that a group extension $E$ of $G$ by $N$ is a short exact
sequence in the category of groups:
{\renewcommand{\theequation}{$\C{E}$}
\begin{equation}\label{sesog}
1\to N\overset{j}{\to} E\overset{p}{\to} G\to 1
\end{equation}
}
It is a pointed fibered object with distinguished fiber $N\cong p^{-1}(1)$.
It is natural to categorify the discrete group $G$ as a base $\C{D}_G$,
using discrete categorification and the group $N$ as a fiber $\C{S}_G$,
using simplicial categorification.

%
% *************** Bundle categorification ************************
%
\subsection{Bundle categorification}
There is a natural lift of the group extension $\C{E}$ to groupoids,
defining the categorification $\C{B}_E$ of $E$ as a disjoint fibration
over $\C{D}_G$ of the simplicial groupoids associated to the fibers of $E$.

The multiplication in $E$ extends in a unique way to a monoidal product.
Similarly, the maps $j$ and $p$ extend in a unique way
to strict monoidal functors.

Then $\C{B}_E\to \C{D}_G$ is a fibered monoidal category,
part of a (short exact) sequence of strict monoidal groupoids
with unit.
Morphisms of group extensions define in an obvious way monoidal functors
between the corresponding fibered categories, and
the correspondence $(E\to G) \mapsto (\C{B}_E\to \C{D}_G)$ is functorial.
\begin{defin}
The functor associating the monoidal groupoid $\C{B}_E$
to a group extension $\C{E}$ is called {\em bundle categorification}
and it is denoted by $\C{B}$.
\end{defin}
The categorification functors are clearly compatible in the
following sense.
\begin{theorem}\label{T:categ}
Discrete, Simplicial and Bundle Categorification maps $\C{D}, \C{S}$ and
$\C{B}$ are functors from the category of groups and
group extensions $\C{E}xt$,
to the category of monoidal groupoids.

$\C{D}$ and $\C{S}$ are the restrictions of $\C{B}$
corresponding to the two natural embeddings:
$$ G \overset{T_b}{\mapsto} (1\to G\to G\to 1\to 1)\quad (trivial\ base)$$
$$ G \overset{T_f}{\mapsto} (1\to 1\to G\to G\to 1)\quad (trivial\ fiber)$$
embedding the category of groups into the category of group extensions:
$$
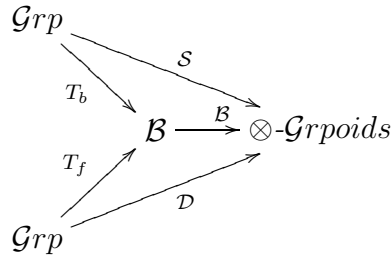
\diagram
\C{G}rp \drto_{T_b} \drrto^{\C{S}} \\
& \C{B} \rto^{\C{B}\quad } & \otimes\text{-}\C{G}rpoids \\
\C{G}rp \urto^{T_f} \urrto_{\C{D}}
\enddiagram$$
\end{theorem}

%
% *************** Relation to group cohomology ************************
%
\subsection{Relation to group cohomology}\label{applcoh}

Consider the group extension \ref{sesog}.
Categorifying the extension as defined in the previous section,
a set-theoretic section is equivalent to
a splitting in the category $\otimes\text{-}\C{G}rpoid$,
as explained below.

Choose such a section $s:G\to E$ of $p$.
There is a unique functor over the section $s$, denoted by $S$.
Since the fiber of $\C{B}_E$ is simplicial, there is a unique
natural isomorphism $\eta$ between the functors $s\circ\otimes$ and
$s^*(\otimes)$:
$$\eta(x,y):S(x\otimes y)\to S(x)\otimes S(y), \quad x,y\in G=Ob(\C{D}_G)$$
Define the functor $\C{R}_E:\C{B}_E \to \C{C}_N$,
the fiberwise analog of the covering transformation
defined in section \ref{covtr}:
$$ \C{R}_E(\overset{x}{\underset{y}{\downarrow}})=n\
\text{iff}\ j(n)=yx^{-1}$$
The natural isomorphism defines what in group language is called
the {\em factor set} of the function $s$ \cite{EML2}:
$$f(x,y)=\C{R}_E(\eta(x,y))=s(x)s(y) s(xy)^{-1}$$
Since any diagram in $\C{B}_E$ commutes,
the functorial morphism $\eta$ is a monoidal structure:
$$\diagram
S((a b)c) \ar@{=}[d] \rto^{\eta(a b,c)} &
S(a b)\otimes S(c) \rto^{\eta(a,b)\otimes I_{S(c)}\quad} &
  (S(a)\otimes S(b))\otimes S(c)\ar@{=}[d]\\
S(a(b c)) \rto^{\eta(a,b c)} & S(a)\otimes S(b c)
\rto^{I_{S(a)}\otimes \eta(b,c)\quad} & S(a)\otimes(S(b)\otimes S(c))
\enddiagram$$
and the pair $(S,\eta)$ is a monoidal functor.
In the above diagram the monoidal product in $G$ was denoted as
concatenation. 

Inner conjugation $C:E\to Aut(E)$ in $E$ defines a left quasi-action
$L:G\to Aut(N)$ of $G$ on $N$: $L(x)(n)=C_{s(x)}(n)$.
Recall that the covering transformation $\C{R}_E$ intertwines
left multiplication by identity morphisms with inner conjugation
and right multiplication by identity morphisms with
the trivial right action of $G$ (see \ref{covtr}).
Applying the functor $\C{R}_E$ to the above monoidality diagram,
one obtains that the factor set $f$ is a group 2-cocycle
relative to the left quasi-action $L$:
$$f(ab,c)f(a,b)=L_a(f(b,c))f(a,bc)$$
Isomorphic group extensions inducing the same quasi-action of $G$ on $N$
define cohomologous 2-cocycles. To see this, 
consider an isomorphism of extensions $\gamma':E\to E'$,
i.e. acting trivial on $N$ and inducing identity on $G$.
It defines a section $s'=\gamma'\circ s$ of $E'$ (functor $S'$),
and a corresponding factor set  $f'$ (monoidal structure $\eta'$).
The induced quasi-action $L'(g)(n)=s'(g) n s'(g)^{-1}$ has the same
conjugacy class $[L']=[L]$ as the one induced by the section $s$.
Since $\gamma$ is a group homomorphism, the unique functor $\Gamma$
extending $\gamma$ on objects is a strict monoidal functor,
and $(S',\eta')=(\Gamma, id) \circ (S,\eta)$ as monoidal functors.
Then the two 2-cocycles $f$ and $f'$
(relative to the different quasi-actions $L$ and $L'$), coincide.

To prove that the cohomology class $[f]$ corresponding to the
extension $E$ is well defined, we may thus assume $E'=E$.
If $s'$ is another section inducing the same outer action $[L]$,
define the 1-cochain $\gamma(x):G\to N$ by:
$$\gamma(x)=n \quad \text{iff}\quad j(n)=s(x) s'(x)^{-1}$$
Then $\gamma$ corresponds through the covering transformation
to the unique natural isomorphism $\Gamma(x):S(x)\to S'(x)$.
Since any diagram in $\C{B}_E$ commutes:
$$\xymatrix @C=2pc @R=.5pc {
S(ab) \ddto_{\eta(ab)} \rto^{\Gamma(ab)} & S'(ab) \ddto^{\eta'(ab)}\\
& & \C{D}_G \ar@/^1pc/[rr]^{(S,\eta)}  \ar@{}[rr]|{\Downarrow \Gamma}
    \ar@/_1pc/[rr]_{(S',\eta')} & & \C{B}_E \\
S(a)\otimes S(b) \rto^{\Gamma(a) \otimes \Gamma(b)} &
  S'(a)\otimes S'(b)
}$$
and $\Gamma$ is an isomorphism of monoidal functors.
The relation obtained by applying the covering transformation:
$$f'(ab)\gamma(ab)=\C{R}_E(\Gamma(a)\otimes \Gamma(b)) f(ab)$$
shows that the two 2-cocycles $(f,L)$ and $(f',L')$ are
{\em weak cohomologous} as defined in \cite{Io}:
$$f'(x,y)\ \partial^-_L\gamma=\partial^+_{L'}\gamma \ f(x,y)$$
If $s$ and  $s'$ induce the same quasi-action $L=L'$, then $\gamma$
is central and one obtains the usual cohomology relation \cite{EML2, Io}:
$$f'(ab)\gamma(ab)=\gamma(a)\gamma(b) f(ab)$$
\begin{rem}
The above categorical interpretation of groups and functions suggests
to consider the full subcategory of group objects in the category {\em Sets}
rather then just groups and group homomorphisms. Then the natural
multiplication of group valued functions provides an internal $Hom$.
\end{rem}

\section{Categorification and topology}

\subsection{Nerve of a category}
The nerve of a category is another example of categorification,
as it will be explain below.

Let $\C{C}$ be a small category.
The categorical analog of the {\em model spaces} of singular
homology are the {\em semi-simplicial} category denoted $\Delta_n$
associated to the total ordered sets $\{0,1,..., n\}$ \cite{Se}.
The analog of singular n-cochains are the $\C{C}$-valued functors
$C^n(\C{C})$ defined on $\Delta_n$.
The semi-simplicial structure coboundary and degeneracy functors
are defined by duality:
$$\partial^i c=c \circ \partial_i, \quad \epsilon^i c=c\circ \epsilon_i$$
with the natural definition of boundary and degeneracy functors defined
on $\Delta_n$.

Let $\Delta$ the {\em quiver of categories and functors}
defined by the above categories $\Delta_n$ and simplicial functors.
The {\em nerve} $N\C{C}$ of the category $\C{C}$ is
the value on $\C{C}$ of the representation functor $Hom(\Delta, \cdot)$
of $\Delta$.
% the sequence of $Hom(\Delta_n,\C{C})$ and maps...

Thus the nerve of a category is a categorification of
acyclic model from topology, and it may be characterized as
the ``singular homology'' of the category $\C{C}$''.
The correspondence is given by
the {\em geometric realization} $B\C{C}$ of the semi-simplicial set $N\C{C}$
\cite{Mi}, which is called the {\em classifying space} of $\C{C}$.

As an example,
the nerve of $\C{S}_G$ is the simplicial structure underlying the
bar construction for the group $G$.
Its geometric realization $EG$ is a contractible, free $G$-space \cite{Se}.

The classifying space $B\C{C}_G$ of the tautological groupoid associated
to the group $G$ is a classifying space $BG$ for $G$ as a discrete group.
All its homotopy groups are trivial except
the fundamental group which is $G$.
At the level of groupoids $\C{S}_G/G\cong \C{C}_G$ through the
covering transformation, as well as for the corresponding topological spaces
$EG/G\cong BG$ \cite{Se}, by taking their geometric realization.

% Q: commutative diagram?, is EG the universal bundle over BG?

\subsection{Topological spaces as categories}

The above procedure associates a topological space to a small category.
We will briefly recall the categorification which,
in the other direction, associates a category to a topological space.

Let $(X,\tau)$ be a topological space. Define the category $\C{C}_X$
of open sets with inclusions as morphisms.
Intersection of open sets is a monoidal product. It is a category
with final object $X$ and with direct limits.
Continuous functions $f$ define functors $f^{-1}$ in the obvious way,
and the construction is functorial.

Families of open sets are the (full) subcategories and covers are
{\em cofinal} subcategories. As usual ``subobjects'' should be
viewed as monomorphisms, so covers will be viewed as
full and faithful functors $U:I\to \C{C}_X$,
defined on posets ($i\le j$ iff $U_i\to U_j$).

Then morphisms of covers:
$$\xymatrix @C=2pc @R=.5pc {
J \ddto_{V} \rto^{\phi} & I \ddto^{U}\\
\ar@{}[r]|{\overset{r}{\Rightarrow}} & &
J \ar@/^1pc/[rr]^{f^*U}  \ar@{}[rr]|{\Downarrow r}
    \ar@/_1pc/[rr]_{V} & & \C{C}_X \\
\C{C}_X \rto^{id} & \C{C}_X
}$$
are precisely {refinements of open covers}.
Alternatively, the map $r$ is a natural transformation between
the pull-back cover $\phi^*U$ and the cover $U$.

The above categorification shows that, in many cases,
2-arrows more general then natural
transformations are needed. The usual approach to $n$-categories,
the ``globular approach'':
$$\diagram
A \ar@/^1pc/[r]^{f}  \ar@{}[r]|{\Downarrow \eta}
    \ar@/_1pc/[r]_{g} &
B \ar@/^1pc/[r]^{f'}  \ar@{}[r]|{\Downarrow \eta'}
    \ar@/_1pc/[r]_{g'} & C
\enddiagram$$
corresponds to interpreting 2-arrows
as {\em fix end homotopies}, while more general transformations
corresponding to arbitrary homotopies are often present:
$$\diagram
\cdot \dto_{\partial^-H} \rrto^{source} &
  \ar@{}[d]|{\Downarrow H} &
  \cdot \dto^{\partial^+H}\\
\cdot \rrto_{target} & & \cdot
\enddiagram$$

%=======================================================================
%                      Bibliography
%=======================================================================

\end{document}